\def\version{0.44}

\def\journal{CA}
\def\titlep{C$^{*}$-bialgebra defined as the direct sum of 
UHF algebras}
\documentclass[11pt]{article}
\usepackage{graphicx,ifthen}
\usepackage{amssymb}
\usepackage{amsmath}

\font\germ=eufm10 at12pt

\def\goth#1{\hbox{\germ#1}}

\newcommand{\qed}{\hbox{\rule[-2pt]{3pt}{6pt}}}
\newcommand{\qedh}{\hfill\qed \\}

\newcommand{\vv}{\vspace{.3in}}

\newcommand{\vep}{\varepsilon}

\def\labelenumi{\theenumi}
\def\theenumi{\arabic{enumi}}
\def\labelenumi{\theenumi}
\def\theenumi{\Alph{enumi}}
\renewcommand{\theenumi}{\Alph{enumi}}
\def\labelenumi{\theenumi}
\def\theenumi{\arabic{enumi}}
\def\labelenumi{\theenumi}
\def\theenumi{{\rm (\roman{enumi})}}
\newtheorem{Thm}{Theorem}[section]

\newtheorem{rem}[Thm]{Remark}

\newtheorem{ex}[Thm]{Example}
\newtheorem{defi}[Thm]{Definition}
\newtheorem{lem}[Thm]{Lemma}

\newtheorem{prop}[Thm]{Proposition}
\newtheorem{prob}[Thm]{Problem}
\newtheorem{cor}[Thm]{Corollary}
\newtheorem{fact}[Thm]{Fact}

\newtheorem{fig}[Thm]{Figure}

\newcommand{\kn}{\Large\bf
$K\hspace{-.4cm} N$
\Large\bf\vv }

\def\cal#1{\mathcal #1}
\def\con{{\cal O}_{n}}
\def\coni{{\cal O}_{\infty}}

\def\pr{{\it Proof.}\quad}

\def\nset#1{\{1,\ldots,n\}^{#1}}

\def\co#1{{\cal O}_{#1}}
\def\disp#1{{\displaystyle #1}}
\def\brl{branching law}
\def\bfsnl{{\rm BFS}_{N}(\Lambda)}

\addtocounter{footnote}{1}
\def\sftt#1{
\setcounter{equation}{0}
\addtocounter{footnote}{1}
\section{#1}
}

\def\ssft#1{\subsection{#1}}

\def\sssft#1{\subsubsection{#1}}

\begin{document}
%
% Personal data
%
\def\autherp{Katsunori Kawamura}
\def\emailp{e-mail: kawamura@kurims.kyoto-u.ac.jp.}
\def\addressp{{\small {\it College of Science and Engineering, 
Ritsumeikan University,}}\\
{\small {\it 1-1-1 Noji Higashi, Kusatsu, Shiga 525-8577, Japan}}
}

\def\infw{\Lambda^{\frac{\infty}{2}}V}
\def\zhalfs{{\bf Z}+\frac{1}{2}}
\def\ems{\emptyset}
\def\pmvac{|{\rm vac}\!\!>\!\! _{\pm}}
\def\vac{|{\rm vac}\rangle _{+}}
\def\dvac{|{\rm vac}\rangle _{-}}
\def\ovac{|0\rangle}
\def\tovac{|\tilde{0}\rangle}
\def\expt#1{\langle #1\rangle}
\def\zph{{\bf Z}_{+/2}}
\def\zmh{{\bf Z}_{-/2}}
\def\brl{branching law}
\def\bfsnl{{\rm BFS}_{N}(\Lambda)}
\def\scm#1{S({\bf C}^{N})^{\otimes #1}}
\def\mqb{\{(M_{i},q_{i},B_{i})\}_{i=1}^{N}}
\def\zhalf{\mbox{${\bf Z}+\frac{1}{2}$}}
\def\zmha{\mbox{${\bf Z}_{\leq 0}-\frac{1}{2}$}}
\newcommand{\mline}{\noindent
\thicklines
\setlength{\unitlength}{.1mm}
\begin{picture}(1000,5)
\put(0,0){\line(1,0){1250}}
\end{picture}
\par
 }
\def\ptimes{\otimes_{\varphi}}
\def\delp{\Delta_{\varphi}}
\def\delps{\Delta_{\varphi^{*}}}
\def\gamp{\Gamma_{\varphi}}
\def\gamps{\Gamma_{\varphi^{*}}}
\def\sem{{\sf M}}
\def\hdelp{\hat{\Delta}_{\varphi}}
\def\tilco#1{\tilde{\co{#1}}}
\def\ndm#1{{\bf M}_{#1}(\{0,1\})}
\def\cdm#1{{\cal M}_{#1}(\{0,1\})}
\def\tndm#1{\tilde{{\bf M}}_{#1}(\{0,1\})}
\def\sck{{\sf CK}_{*}}
\def\hdel{\hat{\Delta}}
% Boldfont
\def\ba{\mbox{\boldmath$a$}}
\def\bb{\mbox{\boldmath$b$}}
\def\bc{\mbox{\boldmath$c$}}
\def\bd{\mbox{\boldmath$d$}}
\def\be{\mbox{\boldmath$e$}}
\def\bk{\mbox{\boldmath$k$}}
\def\bm{\mbox{\boldmath$m$}}
\def\bn{\mbox{\boldmath$n$}}
\def\bp{\mbox{\boldmath$p$}}
\def\bq{\mbox{\boldmath$q$}}
\def\bu{\mbox{\boldmath$u$}}
\def\bv{\mbox{\boldmath$v$}}
\def\bw{\mbox{\boldmath$w$}}
\def\bx{\mbox{\boldmath$x$}}
\def\by{\mbox{\boldmath$y$}}
\def\bz{\mbox{\boldmath$z$}}
\def\bomega{\mbox{\boldmath$\omega$}}
\def\N{{\bf N}}
\def\lxm{L_{2}(X,\mu)}
%%%%%%%%%%%%%%%%%%%%%%%%%%%%%%%%%%%%%%%%%%%%%%%%
\def\titlepage{%\vspace{-4cm}

\noindent
{\bf 
\noindent
\thicklines
\setlength{\unitlength}{.1mm}
\begin{picture}(1000,0)(0,-300)
\put(0,0){\kn \knn\, for \journal\, Ver.\version}
\put(0,-50){\today}
\end{picture}
}
\vspace{-2.3cm}
\quad\\
{\small file: \textsf{tit01.txt,\, J1.tex}
 \footnote{
  ${\displaystyle
  \mbox{directory: \textsf{\fileplace}, 
  file: \textsf{\incfile},\, from \startdate}}$
          }
}
\quad\\
\framebox{
 \begin{tabular}{ll}
 \textsf{Title:} &
 \begin{minipage}[t]{4in}
 \titlep
 \end{minipage}
 \\
 \textsf{Author:} &\autherp
 \end{tabular}
}
%\mline
{\footnotesize	\tableofcontents }
}
%%%%%%%%%%%%%%%%%%%%%%%%%%%%%%%%%%%%%%%%%%%%%%%%
\def\ngt{{\bf N}_{\geq 2}}
\def\ngti{\ngt^{\infty}}
\def\tngti{\widetilde{{\bf N}}_{\geq 2}^{\infty}}
\def\bbn{{\Bbb N}}
\def\bbz{{\Bbb Z}}
\def\bbc{{\Bbb C}}
\def\bbr{{\Bbb R}}
\def\bbf{{\Bbb F}}
\def\bbt{{\Bbb T}}
\def\bbb{{\Bbb B}}

%
%%%%%%%%% Cut from here %%%%%%%%%%
%\input comm.txt
%%%%%%%%% End of Cut %%%%%%%%%
%
%
\setcounter{section}{0}
\setcounter{footnote}{0}
\setcounter{page}{1}
\pagestyle{plain}

%
% Title
%
%%%%%%%%%%%%%%%%%%%%%%%%%%%
\title{\titlep}
\author{\autherp\thanks{\emailp}
\\
\addressp}
\date{}
\maketitle

%%%%%%%%%%%%%%%%%%%%%%%%%%%%%%%%%%%%%%%%
%
% Abstract
%
\begin{abstract}
Let ${\cal A}_{0}(*)$ denote the
direct sum of a certain set of UHF algebras
and let ${\cal A}(*)\equiv {\bf C}\oplus {\cal A}_{0}(*)$.
We introduce a non-cocommutative comultiplication
$\Delta_{\varphi}$ on ${\cal A}(*)$, 
and give an example of comodule-C$^{*}$-algebra 
of the C$^{*}$-bialgebra $({\cal A}(*),\Delta_{\varphi})$.
With respect to $\Delta_{\varphi}$,
we define a non-symmetric tensor product
of $*$-representations of UHF algebras and 
show tensor product formulas of GNS representations 
by product states.
\end{abstract}


%%%%%%%%%%%%%%%%%%%%%%%%%%%%%%%%%%%%%%%%%%%%%%
%
% Section 1
%
\sftt{Introduction}
\label{section:first}
A C$^{*}$-bialgebra is a generalization of bialgebra 
in the theory of C$^{*}$-algebras,
which was introduced in C$^{*}$-algebraic framework 
for quantum groups \cite{KV,MNW}. 
We have studied C$^{*}$-bialgebras and their construction method,
and computed non-symmetric tensor products of $*$-representations 
with respect to non-cocommutative comultiplications 
\cite{TS01,TS02,TS05,TS21,TS24,TS07,TS27}.
In this paper, 
we introduce a non-cocommutative C$^{*}$-bialgebra defined as
the direct sum of a certain set of uniformly hyperfinite (=UHF) algebras.
With respect to the comultiplication,
we define a non-symmetric tensor product
of $*$-representations of UHF algebras, and
show tensor product formulas of 
Gel'fand-Na\u{\i}mark-Segal (=GNS) representations by product states.
The part of tensor product formulas has been given in the previous paper \cite{TS18}
without C$^{*}$-bialgebra.
The present version is reorganized such that 
the tensor product is given by the comultiplication of a C$^{*}$-bialgebra.
In this section, we show our motivation,
definitions and construct the C$^{*}$-bialgebra.
The main feature of this paper is as follows:
\begin{itemize}
\item
A new non-commutative and non-cocommutative C$^{*}$-bialgebra 
is obtained.
In our previous research, we treated only C$^{*}$-bialgebras
$(A,\Delta)$ which satisfy $\Delta(A)\subset A\otimes A$.
In this paper, this property does not holds.
The bialgebra structure does not appear 
unless one takes the direct sum of all UHF algebras.
Until now,
there is no theory which treat all UHF algebras at once.
\item
The C$^*$-bialgebra is naturally constructed by using 
a well-known structure of UHF algebras.
The standard parametrization of
GNS representations of UHF algebras by product states
is compatible with the tensor product formulas.
\item
Tensor product formulas of non-type I representations
are obtained for the first time except \cite{TS07}.
\item
A construction method of 
C$^{*}$-bialgebra is a little bit generalized.
\end{itemize}
%%%%%%%%%%%%%%%%%%%%%%%%%%%%%%%%%%%%%%%%%%%%%%%%%%%%%%%%%%%%%%%%%%%%%%%%
%
% subsection 1.1
%
\ssft{Motivation}
\label{subsection:firstone}
In this subsection, we explain our motivation 
and the background of this study.
Explicit mathematical definitions will 
be shown after $\S$ \ref{subsection:firsttwo}.

According to \cite{FH,Herschend,Pukanszky},
given two representations of a group $G$, 
their tensor product (or Kronecker product
\cite{Pukanszky}) is a new representation
of $G$, which decomposes into a direct sum of indecomposable representations.
The problem of finding this decomposition is called 
the {\it Clebsch-Gordan problem}
and the resulting formula for the decomposition is called 
the {\it tensor product formula}
(or {\it Clebsch-Gordan formula} \cite{Herschend}).
Furthermore, the tensor product
is important to describe the duality of $G$ \cite{Tatsuuma}.
A generalization of the Clebsch-Gordan problem for groups is to consider
modules over associative algebras instead of group algebras. However, there lies an
obvious obstruction in that there is no known way to define the tensor product of
two left modules over an arbitrary associative algebra. 
For group algebras, 
the extra structure coming from the group yields the tensor product.
For a bialgebra $A$, 
the associative tensor product of representations of $A$ can be
defined by using the comultiplication. 
In this way, one of most important motivations of 
the study of bialgebras is the tensor product of their representations.

In \cite{TS01}, we introduced a non-symmetric tensor product 
among all $*$-representations of Cuntz algebras and determined tensor product 
formulas of all permutative representations completely,
in spite of the unknown of any comultiplication of Cuntz algebras.
In \cite{TS02}, we generalized this construction of 
tensor product to a system of C$^{*}$-algebras
and $*$-homomorphisms indexed by a monoid.
For example, we constructed a non-symmetric tensor product
of all $*$-representations of Cuntz-Krieger algebras 
by using Kronecker products of matrices \cite{TS11,TS05,TS07}.

On the other hand,
UHF algebras and their $*$-representations 
are well studied \cite{AK02R,AW,A,AN,BJ,Glimm, Powers,N}.
For example,
GNS representations of product states of UHF algebras 
were completely classified by \cite{AN}.
This class contains $*$-representations of UHF algebras
of all Murray-von Neumann's types I, II, III (\cite{Blackadar2006}, $\S$ III.5).

Our interests are to construct a C$^{*}$-bialgebra from
UHF algebras, and to define a tensor product of $*$-representations of UHF algebras 
with respect to the comultiplication.
Since one knows neither cocommutative nor non-cocommutative comultiplication 
of UHF algebras, the tensor product is new if one can find it.

%%%%%%%%%%%%%%%%%%%%%%%%%%
%
% subsection 1.2
%
\ssft{C$^{*}$-bialgebra}
\label{subsection:firsttwo}
In this subsection,
we recall terminology about C$^{*}$-bialgebra according to \cite{ES,KV,MNW}.
For two C$^{*}$-algebras $A$ and $B$,
we write ${\rm Hom}(A,B)$ as the set of all $*$-homomorphisms 
from $A$ to $B$,
and let ${\cal M}(A)$ denote the multiplier algebra of $A$.
We assume that every tensor product $\otimes$ 
as below means the minimal C$^{*}$-tensor product.
%
% Definition 1.1
%
\begin{defi}
\label{defi:cstar}
A pair $(A,\Delta)$ is a C$^{*}$-bialgebra
if $A$ is a C$^{*}$-algebra and 
$\Delta\in {\rm Hom}(A,{\cal M}(A\otimes A))$
such that the linear span of $\{\Delta(a)(b\otimes c):a,b,c\in A\}$ 
is norm dense in $A\otimes A$ and 
$(\Delta\otimes id)\circ \Delta=(id\otimes\Delta)\circ \Delta$.
We call $\Delta$ the comultiplication of $A$.
\end{defi}

\noindent
We say that a C$^{*}$-bialgebra $(A,\Delta)$ is {\it unital}
if $A$ is unital and $\Delta$ is unital;
$(A,\Delta)$ is {\it counital}
if there exists $\vep\in {\rm Hom}(A,{\bf C})$ which satisfies
$(\vep\otimes id)\circ \Delta= id = (id\otimes \vep)\circ\Delta$.
We call $\vep$ the {\it counit} of $A$ and write $(A,\Delta,\vep)$ 
as the counital C$^{*}$-bialgebra $(A,\Delta)$ with the counit $\vep$.
Remark that Definition \ref{defi:cstar} does not mean 
$\Delta(A)\subset A\otimes A$.
A {\it bialgebra} in the purely algebraic theory \cite{Abe,Kassel} means 
a unital counital bialgebra with the unital counit
with respect to the algebraic tensor product,
which does not need to have an involution.
Hence a C$^{*}$-bialgebra is not a bialgebra in general.

According to \cite{TS02},
we recall several notions of C$^{*}$-bialgebra.
A $*$-homomorphism $f$ from $A$ to ${\cal M}(B)$ is 
{\it nondegenerate} if $f(A)B$ is dense in $B$.
A pair $(B,\Gamma)$ is a {\it right comodule-C$^{*}$-algebra}
of a C$^{*}$-bialgebra $(A,\Delta)$
if $B$ is a C$^{*}$-algebra and $\Gamma$ is a nondegenerate 
$*$-homomorphism from $B$ to ${\cal M}(B\otimes A)$ which satisfies
$(\Gamma\otimes id)\circ \Gamma=
(id\otimes \Delta)\circ \Gamma$
where both $\Gamma\otimes id$ and 
$id\otimes \Delta$ are extended to unital $*$-homomorphisms
from ${\cal M}(B\otimes A)$ to ${\cal M}(B\otimes A\otimes A)$.
In this case, 
the map $\Gamma$ is called the {\it right coaction} of $A$ on $B$.
A C$^{*}$-bialgebra $(A,\Delta)$ is {\it proper}
if $\Delta(a)(I\otimes b),(b\otimes I)\Delta(a)\in A\otimes A$
for any $a,b\in A$ where $I$ denotes the unit of ${\cal M}(A)$.
A proper C$^{*}$-bialgebra $(A,\Delta)$ satisfies the {\it cancellation law} 
if $\Delta(A)(I\otimes A)$ and $\Delta(A)(A\otimes I)$ are dense in $A\otimes A$
where $\Delta(A)(I\otimes A)$ and $\Delta(A)(A\otimes I)$ denote 
the linear spans of sets 
$\{\Delta(a)(I\otimes b):a,b\in A\}$ and $\{\Delta(a)(b\otimes I):a,b\in A\}$,
respectively.

%%%%%%%%%%%%%%%%%%%%%%%%%%%%%%%%%%%%%%%%%%%%%%%%%%%%%%%%%%%%%%%%%%%%%%%%%%%%%%
%
% subsection 1.3
%
\ssft{UHF algebras and $*$-isomorphisms among their tensor products}
\label{subsection:firstthree}
In this subsection, we 
 recall UHF algebras \cite{Glimm} and
introduce a set of $*$-isomorphisms among UHF algebras
and their tensor products.

Let ${\bf N}\equiv \{1,2,3,\ldots\}$, $\ngt\equiv \{2,3,4,\ldots\}$ and 
let $\ngti$ denote the set of all sequences in $\ngt$.
For $n\in {\bf N}$,
let $M_{n}$ denote the (finite-dimensional)
C$^{*}$-algebra of all $n\times n$-complex matrices.
For $\ba= (a_{1},a_{2},\ldots)\in \ngti$,
the sequence $\{M_{a_{n}}\}_{n\geq 1}$ of C$^{*}$-algebras
defines the tensor product 
${\cal A}_{n}(\ba)\equiv \bigotimes_{j=1}^{n}M_{a_{j}}$.
With respect to the embedding 
%
% Equation 1.1
%
\begin{equation}
\label{eqn:psi}
\psi_{\ba}^{(n)}:{\cal A}_{n}(\ba)\ni A\mapsto
 A\otimes I\in {\cal A}_{n}(\ba)\otimes M_{a_{n+1}}
={\cal A}_{n+1}(\ba),
\end{equation}
we regard ${\cal A}_{n}(\ba)$ as a C$^{*}$-subalgebra of 
${\cal A}_{n+1}(\ba)$ and
let 
${\cal A}(\ba)$ denote the inductive limit of the system
$\{({\cal A}_{n}(\ba),\psi_{\ba}^{(n)}):n\geq 1\}$:
%
% Equation 1.2
%
\begin{equation}
\label{eqn:ab}
{\cal A}(\ba)\equiv  \varinjlim ({\cal A}_{n}(\ba),\psi_{\ba}^{(n)}).
\end{equation}
By definition, ${\cal A}(\ba)$ is a UHF algebra of 
Glimm's type $\{a_{1},a_{1}a_{2},a_{1}a_{2}a_{3},\ldots\}$ 
which was classified by \cite{Glimm}.
On the contrary, any UHF algebra is isomorphic to ${\cal A}(\ba)$ 
for some $\ba\in \ngti$.
Hence we call ${\cal A}(\ba)$ a UHF algebra in this paper.

%By using the set,
%we will define a tensor product of representations and that of states
%in $\S$ \ref{subsection:firstfour}.

Let $\{E^{(n)}_{i,j}:i,j=1,\ldots,n\}$ denote the set 
of standard matrix units of $M_{n}$.
For $\ba=(a_{n}),\bb=(b_{n})\in \ngti$,
let $\ba\cdot \bb\equiv (a_{1}b_{1},a_{2}b_{2},\ldots)\in\ngti$.
Then $(\ngti,\,\cdot)$ is a commutative semigroup.
For $\ba,\bb\in\ngti$,
define the (standard) $*$-isomorphism $\varphi_{\ba,\bb}^{(n)}$ 
from ${\cal A}_{n}(\ba\cdot \bb)$
onto ${\cal A}_{n}(\ba)\otimes {\cal A}_{n}(\bb)$ by
%
% Equation 1.3
%
\begin{equation}
\label{eqn:embeddingtwo}
\varphi_{\ba,\bb}^{(n)}(E_{j_{1},k_{1}}^{(a_{1}b_{1})}\otimes \cdots \otimes
E_{j_{n},k_{n}}^{(a_{n}b_{n})})
\equiv 
(E_{j_{1}^{'},k_{1}^{'}}^{(a_{1})}\otimes\cdots\otimes 
E_{j_{n}^{'},k_{n}^{'}}^{(a_{n})})
\otimes 
(E_{j_{1}^{''},k_{1}^{''}}^{(b_{1})}\otimes\cdots\otimes 
E_{j_{n}^{''},k_{n}^{''}}^{(b_{n})})
\end{equation}
for each $j_{i},k_{i}\in \{1,\ldots,a_{i}b_{i}\}$, $i=1,\ldots,n$
where
$j_{1}^{'},\ldots,j_{n}^{'}$,
$k_{1}^{'},\ldots,k_{n}^{'}$,
$j_{1}^{''},\ldots,j_{n}^{''}$,
$k_{1}^{''},\ldots,k_{n}^{''}$
are defined as
$j_{i}=b_{i}(j_{i}^{'}-1)+j_{i}^{''}$ and
$k_{i}=b_{i}(k_{i}^{'}-1)+k_{i}^{''}$,
$j_{i}^{'},k^{'}_{i}\in \{1,\ldots,a_{i}\}$,
$j_{i}^{''},k_{i}^{''}\in \{1,\ldots,b_{i}\}$
for each $i=1,\ldots,n$.
For $\psi_{\ba}^{(n)}$ in (\ref{eqn:psi}), we see that
$(\psi_{\ba}^{(n)}\otimes \psi_{\bb}^{(n)})\circ 
\varphi^{(n)}_{\ba,\bb}
=
\varphi^{(n+1)}_{\ba,\bb}\circ \psi_{\ba\cdot \bb}^{(n)}$
for each $\ba,\bb$ and $n$.
From this,
we can define a unique $*$-isomorphism $\varphi_{\ba,\bb}$ 
from ${\cal A}(\ba\cdot \bb)$
onto ${\cal A}(\ba)\otimes {\cal A}(\bb)$ such that 
%
% Equation 1.4
%
\begin{equation}
\label{eqn:isomorphism}
(\varphi_{\ba,\bb})|_{{\cal A}_{n}(\ba\cdot \bb)}=\varphi_{\ba,\bb}^{(n)}
\quad(n\geq 1)
\end{equation}
where we identify ${\cal A}(\ba)\otimes {\cal A}(\bb)$ 
with the inductive limit of the system
$\{({\cal A}_{n}(\ba)\otimes {\cal A}_{n}(\bb),
\psi^{(n)}_{\ba}\otimes \psi^{(n)}_{\bb}):n\geq 1\}$.

We add the unit ${\bf 1}\equiv (1,1,1,\ldots)$ for the semigroup $\ngti$ and
write 
%
% Equation 1.5
%
\begin{equation}
\label{eqn:tngti}
\tngti\equiv \ngti\cup\{{\bf 1}\},
\end{equation}
which is a subsemigroup of ${\bf N}^{\infty}$.
For convenience,
define the $1$-dimensional C$^{*}$-algebra
%
% Equation 1.6
%
\begin{equation}
\label{eqn:onedimension}
{\cal A}({\bf 1})\equiv {\bf C}.
\end{equation}
Remark that ${\cal A}({\bf 1})$ is not a UHF algebra by definition
\cite{Glimm}.
In addition,
for $\ba\in\ngti$,
define 
$\varphi_{{\bf 1},{\bf 1}}$,
$\varphi_{{\bf 1},\ba}$ and
$\varphi_{\ba,{\bf 1}}$ by
%
% Equation 1.7
%
\begin{equation}
\label{eqn:unitdef}
\varphi_{{\bf 1},{\bf 1}}=id_{{\cal A}({\bf 1})},\quad
\varphi_{{\bf 1},\ba}(x)\equiv 1\otimes x,\quad
\varphi_{\ba,{\bf 1}}(x)\equiv x\otimes 1,\quad(x\in {\cal A}(\ba))
\end{equation}
where we identify ${\cal A}({\bf 1})\otimes {\cal A}({\bf 1})$  
with ${\cal A}({\bf 1})$.

%%%%%%%%%%%%%%%%%%%%%%%%%%%%%%%%%%%%%
%
% Remark 1.2
%
\begin{rem}
\label{rem:one}
{\rm
We consider the meaning of (\ref{eqn:embeddingtwo}).
For two matrices $A\in M_{n}$ and $B\in M_{m}$,
define the matrix $A\boxtimes B\in M_{nm}$
by 
%
% Equation 1.8
%
\begin{equation}
\label{eqn:box}
(A\boxtimes B)_{m(i-1)+i^{'},m(j-1)+j^{'}}\equiv 
A_{i,j}B_{i^{'},j^{'}}
\end{equation}
for $i,j\in \{1,\ldots,n\}$
and $i^{'},j^{'}\in \{1,\ldots,m\}$.
The new matrix $A\boxtimes B$ 
is called the {\it Kronecker product} of $A$ and $B$ \cite{Dief,SS}.
For $\ba,\bb\in \ngti$, we see that
%
% Equation 1.9
%
\begin{equation}
\label{eqn:oneone}
\varphi_{\ba,\bb}^{(1)}(A\boxtimes B)
=A \otimes B\quad(A\in M_{a_{1}},\,B\in M_{b_{1}}).
\end{equation}
Hence $\varphi_{\ba,\bb}^{(1)}$ is the inverse operation of the Kronecker product,
which should be called the {\it Kronecker coproduct}.
}
\end{rem}

%%%%%%%%%%%%%%%%%%%%%%%%%%%%%%%%%%5
%
% subsection 1.4
%
\ssft{Main theorems}
\label{subsection:firstfour}
In this subsection, we show our main theorems.
%%%%%%%%%%%%%%%%%%%%%%%%%%%%%%%%%%%%%%%%%%%%%
%
% subsubsection 1.4.1
%
\sssft{Construction of C$^{*}$-bialgebra}
\label{subsubsection:firstfourone}
In this subsubsection,
we construct a C$^{*}$-bialgebra.
For the uncountable set $\{{\cal A}(\ba):\ba\in \tngti\}$ 
of C$^{*}$-algebras in (\ref{eqn:ab})
and (\ref{eqn:onedimension}),
define the direct sum
%
% Equation 1.10
%
\begin{equation}
\label{eqn:astar}
{\cal A}(*)\equiv \bigoplus\{{\cal A}(\ba):\ba\in\tngti\}.
\end{equation}
For $\ba\in \tngti$,
define 
%
% Equation 1.11
%
\begin{equation}
\label{eqn:nset}
{\cal N}_{\ba}\equiv \{(\bb,\bc)\in \tngti\times \tngti:
\bb\cdot\bc=\ba\},
\end{equation}
and 
%
% Equation 1.12
%
\begin{equation}
\label{eqn:atwo}
{\cal A}^{(2)}(\ba)\equiv \bigoplus\{{\cal A}(\bb)\otimes {\cal A}(\bc)
:(\bb,\bc)\in{\cal N}_{\ba}\}.
\end{equation}
We see that
${\cal A}(*)\otimes {\cal A}(*)=\bigoplus\{{\cal A}^{(2)}(\ba):
\ba\in \tngti\}$.
For $\{\varphi_{\ba,\bb}:\ba,\bb\in \tngti\}$
in (\ref{eqn:isomorphism}),
define the $*$-homomorphism $\delp^{(\ba)}$
from ${\cal A}(\ba)$ to ${\cal M}({\cal A}^{(2)}(\ba))$ by
%
% Equation 1.14
%
\begin{equation}
\label{eqn:delpa}
\delp^{(\ba)}(x)\equiv \prod_{(\bb,\bc)\in{\cal N}_{\ba}}\varphi_{\bb,\bc}(x)
\quad(x\in {\cal A}(\ba))
\end{equation}
where we identify ${\cal M}({\cal A}^{(2)}(\ba))$ with 
the direct product 
$\prod_{(\bb,\bc)\in{\cal N}_{\ba}}{\cal A}(\bb)\otimes {\cal A}(\bc)$.
%Since $\varphi_{\bb,\bc}$ is unital,
%$\delp^{(\ba)}$ is also unital.
Define the $*$-homomorphism
$\delp$ from ${\cal A}(*)$ to
${\cal M}({\cal A}(*)\otimes {\cal A}(*))$ by
%
% Equation 1.15
%
\begin{equation}
\label{eqn:detall}
\delp\equiv \bigoplus \{\delp^{(\ba)}:\ba\in \tngti\}
\end{equation}
where 
we also identify $\oplus\{{\cal M}({\cal A}^{(2)}(\ba)):\ba\in \tngti\}$
with 
a C$^{*}$-subalgebra of 
${\cal M}({\cal A}(*)\otimes {\cal A}(*))
\cong \prod \{{\cal M}({\cal A}^{(2)}(\ba)):\ba\in \tngti\}$.
%
% Theorem 1.3
%
\begin{Thm}
\label{Thm:maintwo}
Let $({\cal A}(*),\delp)$ be as in 
(\ref{eqn:astar}) and  (\ref{eqn:detall}), and
let $\vep$ denote the projection
from ${\cal A}(*)$ onto ${\cal A}({\bf 1})$.
Then the following holds:
\begin{enumerate}
%(i)
\item
$({\cal A}(*),\delp)$ is a non-cocommutative 
proper C$^{*}$-bialgebra with counit $\vep$.
%(ii)
\item
$({\cal A}(*),\delp)$ satisfies the cancellation law.
\end{enumerate}
\end{Thm}

%
% Remark 1.5
%
\begin{rem}
\label{rem:subbialgebra}
{\rm
\begin{enumerate}
%(i)
\item
The idea of the definition of $\{\varphi_{\ba,\bb}\}$
in (\ref{eqn:isomorphism})
is an analogy of the set of embeddings of Cuntz algebras
in $\S$ 1.2 of \cite{TS01}.
In $\S$ 6.1 of \cite{TS02},
we also defined a C$^{*}$-bialgebra defined as the direct sum
of a countably infinite set of UHF algebras:
%
% Equation 1.15
%
\begin{equation}
\label{eqn:uhfs}
UHF_{*}\equiv {\bf C}\oplus UHF_2\oplus
UHF_3\oplus UHF_4\oplus\cdots
\end{equation}
where
$UHF_{n}$ is defined as the fixed point subalgebra of 
the Cuntz algebra $\con$ with respect to the $U(1)$-gauge action,
which is $*$-isomorphic onto
the inductive limit $\varinjlim_{k}\, M_{n}^{\otimes k}$.
Clearly, $UHF_{*}$ is a C$^{*}$-subalgebra
of ${\cal A}(*)$ in (\ref{eqn:astar}), but {\it not} 
a C$^{*}$-subbialgebra. The reason is as follows:
The comultiplication $\Delta$ of $UHF_{*}$ in \cite{TS02} satisfies
$\Delta(UHF_{*})\subset UHF_{*}\otimes UHF_{*}$.
On  the other hand,
the restriction $\delp|_{UHF_{*}}$
of the comultiplication $\delp$ of ${\cal A}(*)$ in (\ref{eqn:detall})
satisfies
$\delp(UHF_{*})\not \subset UHF_{*}\otimes UHF_*$.
%%%%%%%%%%%%%%%%%%%%%%%%%%%%%%%%%%%%%%%%%%%%%%%%%%%%%%%%
%(ii)
\item
In order to help reader's understanding,
we demonstrate the image of comultiplication $\delp$
according to the definition.
Let ${\bf 6}=(6,6,\ldots)\in \tngti$.
Then
%
% Equation 1.16
%
\begin{equation}
\label{eqn:nsix}
{\cal N}_{{\bf 6}}=\left\{
({\bf 6},{\bf 1}),
({\bf 1},{\bf 6}),(\ba,\bar{\ba}):
\begin{array}{l}
\ba=(a_{1},a_{2},\ldots),\,
a_{i}\in\{2,3\},\,i\geq 1\\
\bar{\ba}=(6/a_1,6/a_{2},\ldots)
\end{array}
\right\}.
\end{equation}
Remark that ${\cal N}_{{\bf 6}}$
is also a uncountable set.
When $(\bb,\bc)\in {\cal N}_{{\bf 6}}$ and
$\bb=(b_{1},b_{2},b_{3},\ldots)$,
$b_{1}$ is $1$ or $2$ or $3$ or $6$.
Recall 
${\cal A}({\bf 6})=\varinjlim {\cal A}_{n}({\bf 6})
%=UHF_{6}
$
and
${\cal A}_{1}({\bf 6})=M_{6}\supset \{E^{(6)}_{i,j}:i,j=1,\ldots,6\}$.
For $x\in {\cal A}_{1}({\bf 6})$,
%
% Equation 1.17
%
\begin{equation}
\label{eqn:delpx}
\delp(x)=
\delp^{({\bf 6})}(x)
=
\prod_{(\bb,\bc)\in {\cal N}_{{\bf 6}}}
\varphi_{\bb,\bc}(x)
=
\prod_{(\bb,\bc)\in {\cal N}_{{\bf 6}}}
\varphi_{\bb,\bc}^{(1)}(x).
\end{equation}
Since 
$\delp(x)\in 
{\cal M}({\cal A}^{(2)}(\ba))=
\prod_{(\bb,\bc)\in {\cal N}_{{\bf 6}}}
{\cal A}(\bb)\otimes {\cal A}(\bc)$,
we write $\delp(x)$ as an element 
%(uncountably infinite) double sequence 
$(x_{\bb,\bc})_{(\bb,\bc)\in {\cal N}_{{\bf 6}}}$
in $\prod_{(\bb,\bc)\in{\cal N}_{{\bf 6}}}
{\cal A}(\bb)\otimes {\cal A}(\bc)$ (see $\S$ \ref{subsection:secondone}).
Especially,
we compute the case $x=E_{2,2}^{(6)}$.
Since
$E_{2,2}^{(6)}=
E_{1,1}^{(2)}\boxtimes E_{2,2}^{(3)}
=E_{1,1}^{(3)}\boxtimes E_{2,2}^{(2)}$
and (\ref{eqn:oneone}),
components of $\delp(x)$ are given as follows: 
%
% Equation 1.18
%
\begin{equation}
\label{eqn:xbbbc}
x_{\bb,\bc}
=\left\{
\begin{array}{ll}
1\otimes E^{(6)}_{2,2}
\quad &(\mbox{when }b_{1}=1),\\
\\
E_{1,1}^{(2)}\otimes E_{2,2}^{(3)}\quad &(\mbox{when }b_{1}=2),\\
\\
E_{1,1}^{(3)}\otimes E_{2,2}^{(2)}\quad &(\mbox{when }b_{1}=3),\\
\\
E^{(6)}_{2,2}\otimes 1
\quad &(\mbox{when }b_{1}=6).\\
\end{array}
\right.
\end{equation}
This shows that
the flip of the element $\delp(E_{2,2}^{(6)})$ 
in ${\cal M}({\cal A}(*)\otimes {\cal A}(*))$
is not equal to $\delp(E_{2,2}^{(6)})$.
Hence $\delp$ is non-cocommutative.
\end{enumerate}
}
\end{rem}

%%%%%%%%%%%%%%%%%%%%%%%%%%%%%%%%%%%%%%%%%%%%%
%
% subsubsection 1.4.2
%
\sssft{Construction of comodule-C$^{*}$-algebra}
\label{subsubsection:firstfourtwo}
Next, we introduce an example of comodule-C$^{*}$-algebra
of $({\cal A}(*),\delp)$.
Let $\coni$ denote the Cuntz algebra 
generated by the canonical generators $\{s_{i}:i\in {\bf N}\}$ \cite{Cuntz}.
For two sequences $J=(j_{1},\ldots,j_{n}),K=(k_{1},\ldots,k_{n})\in
{\bf N}^{n}$,
define $E_{J,K}^{(\infty)}\equiv s_{j_{1}}\cdots s_{j_{n}}
s_{k_{n}}^{*}\cdots s_{k_{1}}^{*}\in\coni$ 
and define $UHF_{\infty}$ as
the unital C$^{*}$-subalgebra of $\coni$
generated by $\bigcup_{n\geq 1}\{E_{J,K}^{(\infty)}:J,K\in{\bf N}^{n}\}$:
%
% Equation 1.16
%
\begin{equation}
\label{eqn:uhfi}
UHF_{\infty}\equiv C^{*}\langle
\bigcup_{n\geq 1}\{E_{J,K}^{(\infty)}:J,K\in{\bf N}^{n}\}\rangle \subset \coni.
\end{equation}
We prepare some new notations.
For $\ba=(a_{1},a_{2},\ldots)\in \tngti$ and $n\geq 1$,
define the finite subset $S_{n}(\ba)$ of ${\bf N}^{n}$ by
%
% Equation 1.20
%
\begin{equation}
\label{eqn:productset}
%S(\ba)\equiv \bigcup_{n\geq 1}S_n(\ba),\quad
S_{n}(\ba)\equiv \{1,\ldots,a_{1}\}\times \cdots\times
\{1,\ldots,a_{n}\}.
\end{equation}
Let $\{E_{j,k}^{(n)}\}$ be as in $\S$ \ref{subsection:firstthree}.
For $J=(j_{i}),K=(k_{i})\in S_{n}(\ba)$,
define $E_{J,K}^{(\ba)}\in {\cal A}_{n}(\ba)$ by
%
% Equation 1.21
%
\begin{equation}
\label{eqn:ejka}
E_{J,K}^{(\ba)}\equiv E_{j_{1},k_{1}}^{(a_{1})}\otimes \cdots
\otimes E_{j_{n},k_{n}}^{(a_{n})}.
\end{equation}
Define the $*$-homomorphism $\varphi_{\infty,\ba}$
from $UHF_{\infty}$ to $UHF_{\infty}\otimes {\cal A}(\ba)$
by
%
% Equation 1.22
%
\begin{equation}
\label{eqn:ejkd}
\varphi_{\infty,\ba}(E_{J,K}^{(\infty)})
\equiv E_{J^{'},K^{'}}^{(\infty)}\otimes E_{J^{''},K^{''}}^{(\ba)}
\quad(J,K\in {\bf N}^{n},\,n\geq 1)
\end{equation}
where
$J^{'}=(j_{i}^{'}),K^{'}=(k_{i}^{'})\in{\bf N}^{n}$ and $J^{''}=(j_{i}^{''}),
K^{''}=(k_{i}^{''})\in S_{n}(\ba)$
are defined as
$j_{i}=a_{i}(j_{i}^{'}-1)+j_{i}^{''}$ and
$k_{i}=a_{i}(k_{i}^{'}-1)+k_{i}^{''}$ for $i=1,\ldots,n$.
%
% Theorem 1.5
%
\begin{Thm}
\label{Thm:comodulemain}
Let $UHF_{\infty}$ and
$\{\varphi_{\infty,\ba}:\ba\in\tngti\}$ be as in
(\ref{eqn:uhfi}) and (\ref{eqn:ejkd}), respectively.
Define the $*$-homomorphism $\gamp$ from $UHF_{\infty}$
to ${\cal M}(UHF_{\infty}\otimes {\cal A}(*))$ by
%
% Equation 1.23
%
\begin{equation}
\label{eqn:gampprod}
\gamp(x)\equiv \prod_{\ba\in\tngti}\varphi_{\infty,\ba}(x)\quad(x\in UHF_{\infty})
\end{equation}
where we identify
${\cal M}(UHF_{\infty}\otimes {\cal A}(*))$ 
with $\prod \{UHF_{\infty}\otimes {\cal A}(\ba):\ba\in\tngti\}$.
Then $(UHF_{\infty},\Gamma_{\varphi})$
is a right comodule-C$^{*}$-algebra of $({\cal A}(*),\delp)$.
\end{Thm}

In $\S$ \ref{section:second},
we will prove theorems in $\S$ \ref{subsection:firstfour}.
In $\S$ \ref{section:third}, we will show tensor product formulas
among $*$-representations of UHF algebras,
and show more concrete tensor product formulas
for special UHF algebras.

%%%%%%%%%%%%%%%%%%%%%%%%%%%%%%%%%%%%%%%%%%%%%%%%%%%%%%%
%
% Section 2
%
\sftt{Proofs of main theorems}
\label{section:second}
In this section, we show a general method to construct 
C$^{*}$-bialgebras and prove theorems in $\S$ \ref{subsection:firstfour}.
%%%%%%%%%%%%%%%%%%%%%%%%%%%%%%%%%%%%%%%%%%%%%%%%%%%%%%%%%%%%%%%%%%%%%%%
%
% Subsection 2.1
%
\ssft{C$^{*}$-weakly coassociative system}
\label{subsection:secondone}
In this subsection, we review a general method to construct 
C$^{*}$-bialgebras \cite{TS02,TS27},
and generalize it in order to prove theorems.

First, we recall basic facts of
the direct product and the direct sum of 
general C$^{*}$-algebras \cite{Blackadar2006}.
We define two C$^{*}$-algebras $\prod_{i\in\Omega} A_{i}$
and $\bigoplus_{i\in\Omega} A_{i}$ as follows:
$\prod_{i\in\Omega} A_{i}\equiv 
\{(a_{i}):\|(a_{i})\|\equiv \sup_{i}\|a_{i}\|<\infty\}$,
$\bigoplus_{i\in\Omega} A_{i}\equiv 
\{(a_{i}):\|a_{i}\|\to 0\mbox{ as }i\to\infty\}$.
We call $\prod_{i\in\Omega} A_{i}$ and $\bigoplus_{i\in\Omega} A_{i}$ 
the {\it direct product} and the {\it direct sum} of $A_{i}$'s,
respectively. 
It is known that
${\cal M}(\oplus_{i\in \Omega}A_{i})\cong \prod_{i\in \Omega}{\cal M}(A_{i})$.
If $A_{i}$ is unital for each $i$, then
${\cal M}(\oplus_{i\in \Omega}A_{i})\cong \prod_{i\in \Omega}A_{i}$.

Let $\{B_{i}:i\in \Omega\}$ be another set of C$^{*}$-algebras
and let $\{f_{i}:i\in \Omega\}$ be a set of $*$-homomorphisms such that
$f_{i}\in {\rm Hom}(A_{i},B_{i})$ for each $i\in \Omega$.
Then we obtain 
$\oplus_{i\in\Omega} f_{i}\in
{\rm Hom}(\oplus_{i\in\Omega} A_{i},\oplus_{i\in\Omega} B_{i})$.
If $f_{i}$ is nondegenerate for each $i$,
then $\oplus_{i\in \Omega}f_{i}$ is also nondegenerate.
If both $A_{i}$ and $B_{i}$ 
are unital and $f_{i}$ is unital for each $i\in \Omega$,
then $\oplus_{i\in \Omega}f_{i}$ is nondegenerate.

A {\it monoid} is a set $\sem$ equipped with a binary associative operation 
$\sem\times \sem\ni(a,b)\mapsto ab\in \sem$,
and a unit with respect to the operation.
We recall the definition of 
C$^{*}$-weakly coassociative system in \cite{TS27}.
%
% Definition 2.1
% 
\begin{defi}
\label{defi:axiom}
Let $\sem$ be a monoid with the unit $e$.
A data $\{(A_{a},\varphi_{a,b}):a,b\in \sem\}$
is a C$^{*}$-weakly coassociative system (= C$^{*}$-WCS) over $\sem$ if 
$A_{a}$ is a unital C$^{*}$-algebra for $a\in \sem$
and $\varphi_{a,b}$ is a unital $*$-homomorphism
from $A_{ab}$ to $A_{a}\otimes A_{b}$
for $a,b\in \sem$ such that
\begin{enumerate}
%(i)
\item
for all $a,b,c\in \sem$, the following holds:
%
% Equation 2.1
%
\begin{equation}
\label{eqn:wcs}
(id_{a}\otimes \varphi_{b,c})\circ \varphi_{a,bc}
=(\varphi_{a,b}\otimes id_{c})\circ \varphi_{ab,c}
\end{equation}
where $id_{x}$ denotes the identity map on $A_{x}$ for $x=a,c$,
%(ii)
\item
there exists a counit $\vep_{e}$ of $A_{e}$ 
such that $(A_{e},\varphi_{e,e},\vep_{e})$ 
is a counital C$^{*}$-bialgebra,
%(iii)
\item
for each $a\in\sem$, the following holds:
%
% Equation 2.2
%
\begin{equation}
\label{eqn:new}
(\vep_{e}\otimes id_{a})\circ \varphi_{e,a}=id_{a}
=(id_{a}\otimes \vep_{e})\circ \varphi_{a,e}.
\end{equation}
\end{enumerate}
\end{defi}

We slightly generalize Theorem 2.2 in \cite{TS27} as follows.
%
% Theorem 2.2
% 
\begin{Thm}
\label{Thm:mainthree}
Let $\{(A_{a},\varphi_{a,b}):a,b\in \sem\}$ be a C$^{*}$-WCS 
over a monoid $\sem$.
Define C$^{*}$-algebras 
%
% Equation 2.3
%
\begin{equation}
\label{eqn:aca}
A_{*}\equiv  \oplus \{A_{a}:a\in \sem\},\quad
C_{a}\equiv 
\oplus \{A_{b}\otimes A_{c}:(b,c)\in {\cal N}_{a}\}
\quad (a\in\sem) 
\end{equation}
where ${\cal N}_{a}\equiv\{(b,c)\in \sem\times \sem:\,bc=a\}$.
Define $\Delta^{(a)}_{\varphi}
\in{\rm Hom}(A_{a},{\cal M}(C_{a}))$,
$\Delta_{\varphi}
\in {\rm Hom}(A_{*}, {\cal M}(A_{*}\otimes A_{*}))$ and
$\vep\in {\rm Hom}(A_{*},{\bf C})$ by 
%
% Equation 2.4
%
\begin{equation}
\label{eqn:projection}
\left\{
\begin{array}{rl}
\Delta^{(a)}_{\varphi}(x)\equiv &
\disp{
\prod_{(b,c)\in {\cal N}_{a}}
\varphi_{b,c}(x)\quad(x\in A_{a}),}\\
\\
\Delta_{\varphi}\equiv &
\disp{\oplus\{\Delta_{\varphi}^{(a)}:a\in \sem\},}\\
\\
\vep\equiv &\vep_{e}\circ E_{e}
\end{array}
\right.
\end{equation}
where 
we identify
${\cal M}(A_{*}\otimes A_*)$ with $\prod\{{\cal M}(C_{a}):a\in\sem\}$, and
$E_{e}$ denotes the projection from $A_{*}$
onto $A_{e}$.
Then the following holds:
\begin{enumerate}
%(i)
\item
$(A_{*},\delp,\vep)$ 
is a proper counital C$^{*}$-bialgebra.
%(ii)
\item
In addition, if 
$\#{\cal N}_{a}<\infty$ for each $a\in\sem$,
then $\delp(A_{*})\subset A_{*}\otimes A_{*}$
where we naturally identify $A_{*}\otimes A_{*}$
with a C$^{*}$-subalgebra of ${\cal M}(A_{*}\otimes A_{*})$.
\end{enumerate}
\end{Thm}
%
% Proof
%
\pr
(i)
From (\ref{eqn:aca}),
${\cal M}(C_{a})
=\prod \{A_{b}\otimes A_{c}:(b,c)\in {\cal N}_{a}\}$.
Hence $\delp^{(a)}$ is well-defined.
Since 
${\cal M}(A_{*}\otimes A_{*})
=\prod_{a,b\in\sem}A_{a}\otimes A_{b}$,
$\delp$ is also well-defined.
We show the coassociativity of $\delp$.
Let $a,b,c\in\sem$ and
let $Y\in A_{a}\otimes A_{b}\otimes A_{c}$ and $x\in A_{abc}$.
From (\ref{eqn:wcs}) and (\ref{eqn:projection}),
we can verify that
%
% Equation 2.5
%
\begin{equation}
\label{eqn:coaso}
\{(\delp\otimes id)\circ \delp\}(x)Y
=
\{(id\otimes \delp)\circ \delp\}(x)Y.
\end{equation}
This implies 
$\{(\delp\otimes id)\circ \delp\}(x)=\{(id\otimes \delp)\circ \delp\}(x)$
on $A_{*}\otimes A_{*}\otimes A_{*}$
for each $x\in A_{*}$.
Hence the coassociativity is verified.
As the same token,
we can verify that $\vep$ is a counit and 
$(A_*,\delp)$ is proper.

\noindent
(ii)
This follows from Theorem 2.2 of \cite{TS27}.
\qedh

\noindent
We call $(A_{*},\Delta_{\varphi},\vep)$ in 
Theorem \ref{Thm:mainthree} by a (counital)
{\it C$^{*}$-bialgebra} associated with 
$\{(A_{a},\varphi_{a,b}):a,b\in \sem\}$.

We prepare lemmas as follows.
%
% Lemma 2.3
%
\begin{lem}
\label{lem:unitization}
\begin{enumerate}
%(i)
\item
Let $\{(A_{a},\varphi_{a,b}):a,b\in \sem\}$ be a 
C$^{*}$-WCS over a monoid $\sem$
and let $(A_{*},\delp)$ be as in Theorem \ref{Thm:mainthree}
associated with $\{(A_{a},\varphi_{a,b}):a,b\in \sem\}$. 
For $a,b\in\sem$, let $I_{a}$ denote the unit of $A_{a}$ for $a\in\sem$, and 
define
%
% Equation 2.6
%
\begin{equation}
\label{eqn:notation}
X_{a,b}\equiv \varphi_{a,b}(A_{ab})(A_{a}\otimes I_{b}),\quad
Y_{a,b}\equiv \varphi_{a,b}(A_{ab})(I_{a}\otimes A_{b})
\end{equation}
where
$\varphi_{a,b}(A_{ab})(A_{a}\otimes I_{b})$
and $\varphi_{a,b}(A_{ab})(I_{a}\otimes A_{b})$
mean
the linear spans 
of 
$\{\varphi_{a,b}(x)(y\otimes I_{b}):x\in A_{ab},\,y\in A_{a}\}$
and
$\{\varphi_{a,b}(x)(I_{a}\otimes y):x\in A_{ab},\,y\in A_{b}\}$,
respectively.
If both $X_{a,b}$ and $Y_{a,b}$ 
are dense in $A_{a}\otimes A_{b}$ for each $a,b\in \sem$,
then $(A_{*},\delp)$ satisfies the cancellation law.
%(ii)
\item
For a C$^{*}$-WCS $\{(A_{a},\varphi_{a,b}):a,b\in \sem\}$ 
over a monoid $\sem$,
assume that $B$ is a unital C$^{*}$-algebra and 
a set $\{\varphi_{B,a}:a\in \sem\}$ of unital $*$-homomorphisms
such that $\varphi_{B,a}\in {\rm Hom}(B,B\otimes A_{a})$ for each $a\in \sem$ 
and the following holds:
%
% Equation 2.7
%
\begin{equation}
\label{eqn:comodulefour}
(\varphi_{B,a}\otimes id_{b})\circ \varphi_{B,b}
=(id_{B}\otimes \varphi_{a,b})\circ \varphi_{B,ab}\quad (a,b\in \sem).
\end{equation}
Then $B$ is a right comodule-C$^{*}$-algebra of 
the C$^{*}$-bialgebra $(A_{*},\delp)$ with the unital coaction
$\gamp\equiv \prod_{a\in\sem}\varphi_{B,a}$.
\end{enumerate}
\end{lem}
%
% Proof
%
\pr
(i)
By definition,
the algebraic direct sum
$Q\equiv \oplus_{alg}\{\delp(A_{c})(A_{a}\otimes I):c,a\in \sem\}$
is dense in $\delp(A_{*})(A_{*}\otimes I)$.
On the other hand,
%
% Equation 2.8
%
\begin{equation}
\label{eqn:caseac}
\delp(A_{c})(A_{a}\otimes I)=
\left
\{
\begin{array}{ll}
\varphi_{a,b}(A_{ab})(A_{a}\otimes I_{b})\quad &(\exists b\in\sem\,
s.t.\,ab=c),\\
\\
0 \quad &\mbox{(otherwise)}.
\end{array}
\right.
\end{equation}
From this,
%
% Equation 2.9
%
\begin{equation}
\label{eqn:qoplus}
Q=\oplus_{alg}\{\varphi_{a,b}(A_{ab})(A_{a}\otimes I_{b}):a,b\in \sem\}
=\oplus_{alg} \{X_{a,b}:a,b\in \sem\}.
\end{equation}
By assumption, 
$\oplus_{alg} \{X_{a,b}:a,b\in \sem\}$ is dense in 
$\oplus_{alg}  \{A_{a}\otimes A_{b}:a,b\in \sem\}$.
Hence $\delp(A_{*})(A_{*}\otimes I)$ is dense in $A_{*}\otimes A_{*}$.
In a similar fashion,
we see that $\delp(A_{*})(I\otimes A_{*})$ is also dense in $A_{*}\otimes A_{*}$.
Hence the statement holds.

\noindent
(ii)
By identifying
${\cal M}(B\otimes A_{*})$ 
with $\prod\{B\otimes A_{a}:a\in \sem\}$,
$\gamp$ is well-defined.
From (\ref{eqn:comodulefour}), we can verify that
$(\Gamma_{\varphi}\otimes id_{A_{*}})\circ \Gamma_{\varphi}
=(id_{B}\otimes \Delta_{\varphi})\circ \Gamma_{\varphi}$.
Hence the statement holds.
\qedh

%%%%%%%%%%%%%%%%%%%%%%%%%%%%%%%%%%%%%%%%%%5
%
% subsection 2.2
%
\ssft{Proofs of theorems}
\label{subsection:secondtwo}
In this subsection, we prove theorems in $\S$ \ref{subsection:firstfour}.
\\

\noindent
{\it Proof of Theorem \ref{Thm:maintwo}.}
(i)
From Theorem \ref{Thm:mainthree}(i),
it is sufficient to show that
$\{({\cal A}(\ba),\varphi_{\ba,\bb}):\ba,\bb\in\tngti\}$
in
(\ref{eqn:ab}) and (\ref{eqn:isomorphism})
is a C$^{*}$-WCS over the monoid $\tngti$.
By definition,  the following holds:
%
% Equation 2.10
%
\begin{equation}
\label{eqn:commutative}
(\varphi_{\ba,\bb}\otimes id_{\bc})\circ \varphi_{\ba\cdot \bb,\,\bc}
=
(id_{\ba}\otimes \varphi_{\bb,\bc})\circ \varphi_{\ba,\,\bb\cdot \bc}
\quad(\ba,\bb,\bc\in \tngti)
\end{equation}
where $id_{\bx}$ denotes the identity map on ${\cal A}(\bx)$
for $\bx=\ba,\bc$.
Equivalently,
the following diagram is commutative:

\noindent
%%%%%%%%%%%%%%%%%%%%%%%%%%%%%%%%%%%%%
%
%
%
\thicklines
\setlength{\unitlength}{.1mm}
\begin{picture}(1000,450)(-30,50)
\put(-60,450){
\begin{minipage}[t]{2in}
\begin{fig}
\label{fig:first}
\end{fig}
\end{minipage}
}
\put(-30,250){${\cal A}(\ba\cdot \bb\cdot \bc)$}
\put(190,290){\vector(3,2){150}}
\put(150,350){$\varphi_{\ba,\bb\cdot \bc}$}
\put(190,230){\vector(3,-2){150}}
\put(150,160){$\varphi_{\ba\cdot \bb,\bc}$}
\put(390,400){${\cal A}(\ba)\otimes {\cal A}(\bb\cdot \bc)$}
\put(390,100){${\cal A}(\ba\cdot \bb)\otimes {\cal A}(\bc)$}
\put(690,390){\vector(3,-2){150}}
\put(760,350){$id_{\ba}\otimes \varphi_{\bb,\bc}$}
\put(690,120){\vector(3,2){150}}
\put(780,160){$\varphi_{\ba,\bb}\otimes id_{\bc}$}
\put(880,250){${\cal A}(\ba)\otimes {\cal A}(\bb)\otimes {\cal A}(\bc)$.}
\end{picture}

\noindent
Especially, $({\cal A}({\bf 1}),\varphi_{{\bf 1},{\bf 1}})$
is a one-dimensional C$^{*}$-bialgebra with counit $id_{{\cal A}({\bf 1})}$.
Hence 
$\{({\cal A}(\ba),\varphi_{\ba,\bb}):\ba,\bb\in\tngti\}$
is a C$^{*}$-WCS.
The non-cocommutativity of $({\cal A}(*),\delp)$ has been
shown in Remark \ref{rem:subbialgebra}(ii).

\noindent
(ii)
For $\ba,\bb\in \tngti$,
let 
$X_{\ba,\bb}\equiv \varphi_{\ba,\bb}({\cal A}(\ba\cdot \bb))
({\cal A}(\ba)\otimes I_{\bb})$ and
$Y_{\ba,\bb}\equiv \varphi_{\ba,\bb}({\cal A}(\ba\cdot \bb))
(I_{\ba}\otimes {\cal A}(\bb))$.
From Lemma \ref{lem:unitization}(i),
it is sufficient to show that 
both $X_{\ba,\bb}$ and $Y_{\ba,\bb}$ are dense in
${\cal A}(\ba)\otimes {\cal A}(\bb)$.

By definition,
$X_{\ba,\bb}$ is the linear span of 
$\{\varphi_{\ba,\bb}(x)
(y\otimes I_{\bb}):x\in {\cal A}(\ba\cdot \bb),\, y\in {\cal A}(\ba)\}$.
Let $S_{n}(\ba)$ and $E_{J,K}^{(\ba)}$ be as in (\ref{eqn:productset})
and (\ref{eqn:ejka}), respectively.
We see that 
${\cal A}(\ba)\otimes {\cal A}(\bb)$
is linearly spanned by the set
%
% Equation 2.11
%
\begin{equation}
\label{eqn:nmgeq}
\bigcup_{n,m\geq 1}\{E_{J^{'},K^{'}}^{(\ba)}\otimes E_{J^{''},K^{''}}^{(\bb)}:
J^{'},K^{'}\in S_{n}(\ba),J^{''},K^{''}\in S_{m}(\bb)\}
\end{equation}
as a Banach space.

For $n,m\geq 1$,
fix $J^{'},K^{'}\in S_{n}(\ba)$ and $J^{''},K^{''}\in S_{m}(\bb)$.

\def\labelenumi{\theenumi}
\def\theenumi{(\alph{enumi})}
\begin{enumerate}
%(a)
\item
If $n=m$,
then we can choose $J$ and $K$ in $S_{n}(\ba\cdot \bb)$ such that
$E^{(\ba\cdot \bb)}_{J,K}=E_{J^{'},K^{'}}^{(\ba)}\boxtimes E_{J^{''},K^{''}}^{(\bb)}$
where $\boxtimes$ means the componentwise Kronecker product.
Since 
$E_{J^{'},K^{'}}^{(\ba)}\otimes E_{J^{''},K^{''}}^{(\bb)}=
\delp(E^{(\ba\cdot \bb)}_{J,K})( E_{J^{'},K^{'}}^{(\ba)}\otimes I_{\bb})$,
we see 
$E_{J^{'},K^{'}}^{(\ba)}\otimes E_{J^{''},K^{''}}^{(\bb)}
\in \delp({\cal A}(\ba\cdot \bb))({\cal A}(\ba)\otimes I_{\bb})$.
%
%(b)
\item
If $n-m=k>0$,
then we can write as
%
% Equation 2.12
%
\begin{equation}
\label{eqn:ejkbbb}
E_{J^{'},K^{'}}^{(\ba)}\otimes E_{J^{''},K^{''}}^{(\bb)}
=
\sum_{L\in S_{m,k}(\bb)}
E_{J^{'},K^{'}}^{(\ba)}\otimes E_{(J^{''}L),(K^{''}L)}^{(\bb)}
\end{equation}
where 
$S_{m,k}(\bb)\equiv \{1,\ldots,b_{m+1}\}\times
\cdots\times \{1,\ldots,b_{m+k}\}$ 
for $\bb=(b_{1},b_2,\ldots)$
and 
$(J^{''}L)$ denotes the concatenation of two sequences
$J^{''}$ and $L$.
The right hand side of (\ref{eqn:ejkbbb})
is contained in
$\delp({\cal A}(\ba\cdot \bb))({\cal A}(\ba)\otimes I_{\bb})$
from (a).
%(c)
\item
If $n-m=-k<0$,
then this follows from (b) by the same token.
\end{enumerate}
\def\labelenumi{\theenumi}
\def\theenumi{{\rm (\roman{enumi})}}
\def\labelenumii{\theenumii}
\def\theenumii{{\rm (\alph{enumii})}}
%%%%%%%%%%%%%%%%%%%%%%%%%%%%%%%%%%%%%
From (a),(b),(c),
$\delp({\cal A}(\ba\cdot \bb))({\cal A}(\ba)\otimes I_{\bb})$
is dense in ${\cal A}(\ba)\otimes {\cal A}(\bb)$.
Just the same,
we see that 
$\delp({\cal A}(\ba\cdot \bb))(I_{\ba}\otimes {\cal A}(\bb))$
is dense in ${\cal A}(\ba)\otimes {\cal A}(\bb)$.
Hence the statement holds.
\qedh

\noindent
{\it Proof of 
Theorem \ref{Thm:comodulemain}.}
Let $\ba,\bb\in \tngti$, 
$Y\in UHF_{\infty}\otimes {\cal A}(\ba)\otimes {\cal A}(\bb)$  and
$x\in UHF_{\infty}$.
By definition,
we see that
%
% Equation 2.13
%
\begin{equation}
\label{eqn:varphiinf}
\{(\varphi_{\infty,\ba}\otimes id_{\bb})\circ \varphi_{\infty,\,\bb}\}(x)Y
=
\{(id_{\infty}\otimes \varphi_{\ba,\bb})\circ \varphi_{\infty,\,\ba\cdot \bb}\}(x)Y
\end{equation}
where $id_{\infty}$ denotes the identity map on $UHF_{\infty}$.
From this and Lemma \ref{lem:unitization}(ii) for $B=UHF_{\infty}$ 
and $\varphi_{B,\ba}\equiv \varphi_{\infty,\ba}$,
the statement holds.
\qedh

%%%%%%%%%%%%%%%%%%%%%%%%%%%%%%%%%%%%%%%%%%%%%%%%%%%%%%
%
% Section 3
%
\sftt{Tensor product formulas}
\label{section:third}
In this section, we show tensor product formulas
of $*$-representations of the C$^{*}$-algebra ${\cal A}(*)$ with respect to
the comultiplication $\delp$ in $\S$ \ref{subsection:firstfour}.
%%%%%%%%%%%%%%%%%%%%%%%%%%%%%%%%%%%%%%%%%%%%
%
% subsection 3.1
%
\ssft{Basic properties}
\label{subsection:thirdone}
In this subsection, 
we introduce a tensor product
of $*$-representations 
and that of states of UHF algebras,
and show its basic properties.
For a C$^{*}$-algebra ${\goth A}$,
let ${\rm Rep}{\goth A}$
and ${\cal S}({\goth A})$ denote
the class of all $*$-representations
and the set of all states of ${\goth A}$, respectively.
By using the set $\{\varphi_{\ba,\bb}:\ba,\bb\in \ngti\}$
in (\ref{eqn:isomorphism}),
define the operation $\ptimes$ from 
${\rm Rep}{\cal A}(\ba)\times {\rm Rep}{\cal A}(\bb)$
to ${\rm Rep}{\cal A}(\ba\cdot \bb)$ by
%
% Equation 3.1
%
\begin{equation}
\label{eqn:ptimes}
\pi_{1}\ptimes \pi_{2}\equiv (\pi_{1}\otimes \pi_{2})\circ \varphi_{\ba,\bb}
\end{equation}
for $(\pi_{1},\pi_2)\in {\rm Rep}{\cal A}(\ba)\times {\rm Rep}{\cal A}(\bb)$.
We see that 
if $\pi_{i}$ and $\pi_{i}^{'}$
are unitarily equivalent for $i=1,2$,
then 
$\pi_{1}\ptimes \pi_{2}$
and 
$\pi_{1}^{'}\ptimes \pi_{2}^{'}$
are also unitarily equivalent.
Furthermore, 
define the operation $\ptimes$ from 
${\cal S}({\cal A}(\ba))\times {\cal S}({\cal A}(\bb))$
to ${\cal S}({\cal A}(\ba\cdot \bb))$ by
%
% Equation 3.2
%
\begin{equation}
\label{eqn:rho}
\rho_{1}\ptimes \rho_{2}\equiv (\rho_{1}\otimes \rho_{2})\circ \varphi_{\ba,\bb}
\end{equation}
for $(\rho_{1},\rho_2)\in 
{\cal S}({\cal A}(\ba))\times {\cal S}({\cal A}(\bb))$.
From (\ref{eqn:commutative}),
we see that
%
% Equation 3.3
%
\begin{equation}
\label{eqn:ptimestwo}
(\pi_{1}\ptimes \pi_{2})\ptimes \pi_{3}=
\pi_{1}\ptimes (\pi_{2}\ptimes \pi_{3}),\quad
(\rho_{1}\ptimes \rho_{2})\ptimes \rho_{3}=
\rho_{1}\ptimes (\rho_{2}\ptimes \rho_{3})
\end{equation}
for each 
$(\pi_{1},\pi_2,\pi_{3})\in {\rm Rep}{\cal A}(\ba)\times {\rm Rep}{\cal A}(\bb)
\times {\rm Rep}{\cal A}(\bc)$
and 
for $(\rho_{1},\rho_2,\rho_{3})\in
{\cal S}({\cal A}(\ba))\times {\cal S}({\cal A}(\bb))
\times {\cal S}({\cal A}(\bc))$
and $\ba,\bb,\bc\in \ngti$.

The following fact 
is a paraphrase of well-known results of tensor products of factors.
%
% Fact 3.1
%
\begin{fact}
\label{fact:types}
Let $\pi_{1}$ and $\pi_{2}$
be $*$-representations of 
${\cal A}(\ba)$
and 
${\cal A}(\bb)$, respectively.
Then the following holds:
\begin{enumerate}
%(i)
\item
If both $\pi_{1}$ and $\pi_{2}$
are factor representations,
then so is $\pi_{1}\ptimes \pi_{2}$.
%(ii)
\item
The type of $\pi_{1}\ptimes \pi_{2}$
coincides with 
that of $\pi_{1}\otimes \pi_{2}$
where the type of a representation $\pi$ of a C$^{*}$-algebra ${\goth A}$
means the type of the von Neumann algebra $\pi({\goth A})^{''}$
(\cite{Blackadar2006}, Theorem III.2.5.27).
%(iii)
\item
If both $\pi_{1}$ and $\pi_{2}$
are irreducible,
then 
so is $\pi_{1}\ptimes \pi_{2}$.
\end{enumerate}
\end{fact}
%
% Proof
%
\pr
By the definition of $\ptimes$,
%
% Equation 2.1
%
\begin{equation}
\label{eqn:productone}
(\pi_{1}\ptimes \pi_{2})({\cal A}(\ba\cdot \bb))
=
(\pi_{1}\otimes \pi_{2})({\cal A}(\ba)\otimes {\cal A}(\bb)).
\end{equation}
\noindent
(i)
Since $\pi_{1}\otimes \pi_{2}$ is also a factor representation,
the statement holds from (\ref{eqn:productone}).

\noindent
(ii)
By definition,
the type of $\pi_{1}\ptimes \pi_{2}$ is that of 
$\{(\pi_{1}\ptimes \pi_{2})({\cal A}(\ba\cdot \bb))\}^{''}$.
From this and (\ref{eqn:productone}), the statement holds.

\noindent
(iii) By assumption, $\pi_{1}\otimes \pi_{2}$ is also irreducible.
From this and (\ref{eqn:productone}), the statement holds.
\qedh

By definition,
the essential part of the tensor product $\ptimes$
is given by the set $\{\varphi_{\ba,\bb}\}$ of isomorphisms
in (\ref{eqn:embeddingtwo}).
This type of tensor product is known yet
in neither operator algebras nor the purely algebraic theory 
of quantum groups \cite{Kassel}.
%
% Remark 3.2
%
\begin{rem}
\label{rem:second}
{\rm
\begin{enumerate}
%(i)
\item
Our terminology ``tensor product of representations" is different
from usual sense \cite{FH}.
Remark that,
for $\pi,\pi^{'}\in {\rm Rep}{\cal A}(\ba)$,
$\pi\ptimes \pi^{'}\not\in {\rm Rep}{\cal A}(\ba)$ but
$\pi\ptimes \pi^{'}\in {\rm Rep}{\cal A}(\ba\cdot \ba)$
because $\ba\cdot \ba\ne \ba$ for any $\ba\in \ngti$.
%(ii)
\item
From Fact \ref{fact:types}(iii),
there is no nontrivial branch of 
the irreducible decomposition of 
the tensor product of any two irreducibles.
In general, such a tensor product of the other algebra is decomposed
into more than one irreducible component.
For example, see Theorem 1.6 of \cite{TS01}.
\end{enumerate}
}
\end{rem}

%%%%%%%%%%%%%%%%%%%%%%%%%%%%%%%%%%%%%%%%%%
%
% subsection 3.2
%
\ssft{GNS representations by product states and
their tensor product formulas}
\label{subsection:thirdtwo}
We recall well-known GNS representations by product states of UHF algebras
%which are given as infinite tensor products of matrix algebras
\cite{AW,A,N}.
Let $M_{n}$ be as in $\S$ \ref{subsection:firstthree} and 
let $M_{n,+,1}$ denote the set of all positive elements in $M_{n}$ 
whose traces are $1$.
Then a linear functional $\omega$ on $M_{n}$ is a state of $M_{n}$ 
if and only if $\omega$ is equal to the state $\omega_{T}$ which is defined as
$\omega_{T}(x)\equiv {\rm tr}(Tx)$ ($x\in M_{n}$)
for some $T\in M_{n,+,1}$
where ${\rm tr}$ denotes the trace of $M_{n}$.
For $\ba=(a_{1},a_{2},\ldots)\in\ngti$,
let ${\cal T}(\ba)\equiv \prod_{n\geq 1}M_{a_{n},+,1}$.
For a sequence ${\bf T}=(T^{(n)})_{n\geq 1}\in {\cal T}(\ba)$,
define the state
$\omega_{{\bf T}}$ of ${\cal A}(\ba)$ by
%
% Equation 3.5
%
\begin{equation}
\label{eqn:tstate}
\omega_{{\bf T}}(E_{j_{1},k_{1}}^{(a_{1})}
\otimes\cdots\otimes E_{j_{n},k_{n}}^{(a_{n})}
)\equiv T_{k_{1},j_1}^{(1)}\cdots T_{k_{n},j_n}^{(n)}
\end{equation}
for each $j_{1},\ldots,j_{n},k_{1},\ldots,k_{n}$ and $n\geq 1$
where $T^{(n)}_{j,k}$'s denote matrix elements of the matrix $T^{(n)}$.
Then $\omega_{{\bf T}}$ coincides with the product state 
$\bigotimes_{n\geq 1}\omega_{T^{(n)}}$.
%
% Theorem 3.3
%
\begin{Thm}
\label{Thm:awone}(\cite{K-R}, Remark 11.4.16)
For each ${\bf T}\in {\cal T}(\ba)$,
the state $\omega_{{\bf T}}$ in (\ref{eqn:tstate})
is a factor state, that is,
if $({\cal H}_{{\bf T}},\pi_{{\bf T}},\Omega_{{\bf T}})$ is 
the GNS triplet of ${\cal A}(\ba)$
by $\omega_{{\bf T}}$,
then $\pi_{{\bf T}}({\cal A}(\ba))^{''}$ is a factor.
\end{Thm}
The factor ${\cal M}_{{\bf T}}\equiv \pi_{{\bf T}}({\cal A}(\ba))^{''}$ 
is called an {\it Araki-Woods factor (or infinite tensor product 
of finite dimensional type {\rm I} (=ITPFI) factor)} \cite{AW, A}.
Properties of ${\cal M}_{{\bf T}}$ 
and $({\cal H}_{{\bf T}},\pi_{{\bf T}},\Omega_{{\bf T}})$ 
are closely studied in \cite{AW,A,Shlyakhtenko} and \cite{AN}, respectively.

Next, we show tensor product formulas of 
$\pi_{{\bf T}}$'s in Theorem \ref{Thm:awone} as follows.
%
% Theorem 3.4
%
\begin{Thm}
\label{Thm:main}
Let  $\ba,\bb\in \ngti$ and 
let $\omega_{{\bf T}}$ be as in (\ref{eqn:tstate})
with the GNS representation $\pi_{{\bf T}}$.
\begin{enumerate}
%(i)
\item
For each ${\bf T}\in {\cal T}(\ba)$ and ${\bf R}\in {\cal T}(\bb)$, 
%
% Equation 3.6
%
\begin{equation}
\label{eqn:omega}
\omega_{{\bf T}}\ptimes \omega_{{\bf R}}=\omega_{{\bf T}\boxtimes {\bf R}}
\end{equation}
where 
${\bf T}\boxtimes {\bf R}\in {\cal T}(\ba\cdot \bb)$ is
defined as
%
% Equation 3.7
%
\begin{equation}
\label{eqn:boxt}
{\bf T}\boxtimes {\bf R}\equiv (T^{(1)}\boxtimes R^{(1)},
T^{(2)}\boxtimes R^{(2)},T^{(3)}\boxtimes R^{(3)},\ldots)
\end{equation}
for ${\bf T}=(T^{(n)})$ and ${\bf R}=(R^{(n)})$.
%(ii)
\item
For each ${\bf T}\in {\cal T}(\ba)$ and ${\bf R}\in {\cal T}(\bb)$, 
$\pi_{{\bf T}}\ptimes \pi_{{\bf R}}$ 
is unitarily equivalent
to $\pi_{{\bf T}\boxtimes {\bf R}}$.
\end{enumerate}
\end{Thm}
%
% Proof
% 
\pr
(i)
By definition,
the statement holds from direct computation.

\noindent
(ii)
Let $\ba\in \ngti$.
For ${\bf T}\in {\cal T}(\ba)$,
let $({\cal H}_{{\bf T}},\pi_{{\bf T}},\Omega_{{\bf T}})$
denote the GNS triplet by the state $\omega_{{\bf T}}$.
Define the GNS map $\Lambda_{{\bf T}}$  \cite{KV,MNW}
from ${\cal A}(\ba)$ to ${\cal H}_{{\bf T}}$ by
$\Lambda_{{\bf T}}(x)\equiv \pi_{{\bf T}}(x)\Omega_{{\bf T}}$
for $x\in {\cal A}({\bf T})$.
Let $\ba,\bb\in \ngti$.
For ${\bf T}\in {\cal T}(\ba)$ and ${\bf R}\in {\cal T}(\bb)$,
define the unitary $U^{({\bf T},{\bf R})}$
from ${\cal H}_{{\bf T}\boxtimes {\bf R}}$
to ${\cal H}_{{\bf T}}\otimes {\cal H}_{{\bf R}}$ by
%
% Equation 3.8
%
\begin{equation}
\label{eqn:utr}
U^{({\bf T},{\bf R})}
\Lambda_{{\bf T}\boxtimes {\bf R}}(x)
\equiv 
(\Lambda_{{\bf T}}
\otimes
\Lambda_{{\bf R}})(\varphi_{\ba,\bb}(x))
\quad(x\in {\cal A}(\ba\cdot \bb)).
\end{equation}
Since $\varphi_{\ba,\bb}$ is bijective,
$U^{({\bf T},{\bf R})}$ is well-defined as a unitary,
and we see that
%
% Equation 3.9
%
\begin{equation}
\label{eqn:ur}
U^{({\bf T},{\bf R})}\pi_{{\bf T}\boxtimes {\bf R}}(x)
(U^{({\bf T},{\bf R})})^{*}
=
(\pi_{{\bf T}}\ptimes \pi_{{\bf R}})(x)
\quad(x\in {\cal A}(\ba\cdot \bb)).
\end{equation}
Hence two representations
$\pi_{{\bf T}\boxtimes {\bf R}}$ and 
$\pi_{{\bf T}}\ptimes \pi_{{\bf R}}$
are unitarily equivalent.
\qedh

\noindent
From Theorem \ref{Thm:main},
the tensor product $\ptimes$ is compatible with 
product states and their GNS representations.
More precisely, for the following 
two semigroups $({\cal T},\boxtimes)$ and $({\cal S},\ptimes)$, the map 
%
% Equation 3.10
%
\begin{equation}
\label{eqn:calt}
{\cal T}\equiv \bigcup_{\ba\in \ngti}{\cal T}(\ba)
\ni {\bf T}\mapsto \omega_{{\bf T}}\in {\cal S}\equiv 
\bigcup_{\ba\in \ngti}{\cal S}({\cal A}({\ba}))
\end{equation}
is a semigroup homomorphism.
Let ${\sf R}_{\ba}$ denote the set of all unitary equivalence classes
in ${\rm Rep}{\cal A}(\ba)$.
Then 
%
% Equation 3.11
%
\begin{equation}
\label{eqn:caltt}
{\cal T}\ni {\bf T}\mapsto [\pi_{{\bf T}}]\in {\sf R}\equiv 
\bigcup_{\ba\in \ngti}{\sf R}_{{\ba}}
\end{equation}
is also a  semigroup homomorphism from
$({\cal T},\boxtimes)$ to $({\sf R},\ptimes)$
where $[\pi]$ denotes the unitary equivalence class of 
a representation $\pi$.

%%%%%%%%%%%%%%%%%%%%%%%%%%%%%%%%%%%%%%%%%%%%%%%%%%%%%%%%%%%%%%%%%%%%%%%%%%%
%
% Subsection 3.3
%
\ssft{Examples}
\label{subsection:thirdthree}
In this subsection,
we show examples of Theorem \ref{Thm:main}
for special UHF algebras.
Let ${\cal A}(\ba)$ and ${\cal T}(\ba)$ be as in $\S$ \ref{subsection:thirdtwo}.
For $n\geq 1$,
let 
%
% Equation 3.12
%
\begin{equation}
\label{eqn:nnn}
\bn\equiv (n,n,n,\ldots)\in {\bf N}^{\infty}
\end{equation}
and let 
%
% Equation 3.13
%
\begin{equation}
\label{eqn:uhfn}
UHF_{n}\equiv {\cal A}(\bn)=(M_{n})^{\otimes \infty}\quad(n\geq 2).
\end{equation}
Then $UHF_{n}$ is 
the UHF algebra of Glimm's type $\{n^{l}\}_{l\geq 1}$.

Let $\{E^{(n)}_{i,j}\}$ be as in $\S$ \ref{subsection:firstthree}.
For $j\in \{1,\ldots,n\}$,
define $F_{j}^{(n)}\equiv E_{j,j}^{(n)}$.
Then $F_{j}^{(n)}\in M_{n,+,1}$ for each $j$.
For $J=(j_{1},j_2,\ldots)\in \{1,\ldots,n\}^{\infty}$,
define 
${\bf T}(J)\equiv (F_{j_{1}}^{(n)},F_{j_{2}}^{(n)},\ldots)\in {\cal T}(\bn)$.
From (\ref{eqn:tstate}), we see that
%
% Equation 3.14
%
\begin{equation}
\label{eqn:tstatep}
\omega_{{\bf T}(J)}(E_{l_{1},k_{1}}^{(n)}
\otimes\cdots\otimes E_{l_{m},k_{m}}^{(n)}
)\equiv 
\delta_{l_{1},j_1}\cdots \delta_{l_{m},j_m}\,
\delta_{k_{1},j_1}\cdots \delta_{k_{m},j_m}
\end{equation}
for each $l_{1},\ldots,l_{m},k_{1},\ldots,k_{m}\in\{1,\ldots,n\}$ and $m\geq 1$.
For $J=(j_{n})_{n\in {\bf N}},J^{'}=(j^{'}_{n})_{n\in {\bf N}}\in\nset{\infty}$,
we write
$J\approx J^{'}$ if there exists an integer $n_{0}\geq 1$ 
such that $j_{r}=j^{'}_{r}$ for each $r\geq n_{0}$.

%
% Proposition 3.5
%
\begin{prop}
\label{prop:subclass}
For $J\in \{1,\ldots,n\}^{\infty}$,
let $\pi_{{\bf T}(J)}$ denote
the GNS representation of $UHF_{n}$ by the state
$\omega_{{\bf T}(J)}$ in (\ref{eqn:tstatep}), and 
let $P_{n}[J]$ denote  the unitary equivalence class of $\pi_{{\bf T}(J)}$.
Then the following holds:
\begin{enumerate}
%(i)
\item
For each $J\in \{1,\ldots,n\}^{\infty}$,
$\pi_{{\bf T}(J)}$ is irreducible.
%(ii)
\item
For $J,J^{'}\in\nset{\infty}$, $P_n[J]= P_n[J^{'}]$ if and only 
if $J\approx J^{'}$.
%(iii)
\item
For each 
$J=(j_{l})\in \{1,\ldots,n\}^{\infty}$,
$K=(k_{l})\in \{1,\ldots,m\}^{\infty}$
and $n,m\geq 2$, 
%
% Equation 3.15
%
\begin{equation}
\label{eqn:productb}
P_{n}[J]\ptimes P_{m}[K]=P_{nm}[J\star K]
\end{equation}
where
$J\star K\in {\cal T}(\bn\cdot \bm)$ is defined by
%
% Equation 3.16
%
\begin{equation}
\label{eqn:dotdot}
J\star K\equiv (m(j_{1}-1)+k_{1},m(j_{2}-1)+k_{2},m(j_{3}-1)+k_{3},\ldots).
\end{equation}
\end{enumerate}
\end{prop}
%
% Proof
%
\pr
(i)
Since the state $\omega_{{\bf T}(J)}$ is the product state
of pure states,
$\omega_{{\bf T}(J)}$ is also pure.
Hence the statement holds.

\noindent
(ii)
From (2.1) in \cite{AK02R} (see also \cite{BC}),  
the statement holds.

\noindent
(iii)
Recall $\boxtimes$ in (\ref{eqn:boxt}). 
Then we can verify that
${\bf T}(J)\boxtimes {\bf T}(K)={\bf T}(J\star K)$
for each
$J\in \{1,\ldots,n\}^{\infty}$ and
$K\in \{1,\ldots,m\}^{\infty}$.
From this and Theorem \ref{Thm:main}(ii),
the statement holds.
\qedh

\noindent
In Theorem 1.6 of \cite{TS01},
``$J\star K$" is written as the different notation ``$J\cdot K$."

%
% Example 3.6
%
\begin{ex}
\label{ex:noncoco}
{\rm
From Proposition \ref{prop:subclass},
${\sf P}\equiv \bigcup_{n\geq 2}\{P_{n}[J]:J\in\{1,\ldots,n\}^{\infty}\}$
is a semigroup 
of unitary equivalence classes of irreducible representations
with respect to the product $\ptimes$. 
For $n\geq 1$, let $\bn$ be as in (\ref{eqn:nnn}).
Then
%
% Equation 3.17
%
\begin{equation}
\label{eqn:pppp}
P_{2}[{\bf 1}]\ptimes P_{2}[{\bf 2}]=P_{4}[{\bf 2}],\quad 
P_{2}[{\bf 2}]\ptimes P_{2}[{\bf 1}]=P_{4}[{\bf 3}]
\end{equation}
because ${\bf 1}\star {\bf 2}={\bf 2}$
and ${\bf 2}\star {\bf 1}={\bf 3}$.
Since ${\bf 2}\not \approx {\bf 3}$,
$P_{4}[{\bf 2}]\ne P_{4}[{\bf 3}]$.
Hence $({\sf P},\ptimes)$ is non-commutative.
}
\end{ex}

The class $P_{n}[J]$ in Proposition \ref{prop:subclass}(ii)
coincides with
the restriction of a permutative representation of the Cuntz algebra $\con$
on the UHF subalgebra of $U(1)$-gauge invariant elements in $\con$, 
which is called an atom \cite{AK02R,BJ}.
This class contains only type I representations of $UHF_{n}$.
Relations with representations of Cuntz algebras 
and quantum field theory are well studied \cite{AK02R,BJ}.

From Example \ref{ex:noncoco}
and pure states associated with
$P_{2}[{\bf 1}]$ and  $P_{2}[{\bf 2}]$,  
the following holds.
%
% Corollary 3.7
%
\begin{cor}
\label{cor:nonsymmetric}
Let ${\cal S}({\cal A}(\ba))$ and ${\sf R}_{\ba}$
be as in $\S$ \ref{subsection:thirdone}.
\begin{enumerate}
%(i)
\item
There exist 
$\ba,\bb\in\ngti$, and states
$\omega\in {\cal S}({\cal A}(\ba))$ and
$\omega^{'}\in {\cal S}({\cal A}(\bb))$ such that
$\omega\ptimes \omega^{'}\ne \omega^{'}\ptimes \omega$.
%(ii)
\item
There exist 
$\ba,\bb\in\ngti$, and
classes
$[\pi]\in {\sf R}_{\ba}$ and
$[\pi^{'}]\in {\sf R}_{\bb}$
such that
$[\pi]\ptimes [\pi^{'}]\ne 
[\pi^{'}]\ptimes [\pi]$.
\end{enumerate}
\end{cor}

\noindent
From Corollary \ref{cor:nonsymmetric}(i) and (ii),
we say that $\ptimes$ is non-symmetric.
In other words,
two semigroups $({\cal S},\ptimes)$
and $({\sf R},\ptimes)$ are non-commutative.
These non-commutativities come from the
non-commutativity of the Kronecker product of matrices
in Remark \ref{rem:one}.

%
% Problem 3.8
%
\begin{prob}
\label{prob:af}
{\rm
\begin{enumerate}
%(i)
\item
Generalize the tensor product in (\ref{eqn:ptimes})
to that of representations of approximately finite dimensional (=AF) algebras
which are not always UHF algebras.
%(ii)
\item
Reconstruct UHF algebras from 
$({\cal S},\ptimes)$ and $({\sf R},\ptimes)$,
and show a Tatsuuma duality type theorem 
for UHF algebras \cite{Tatsuuma}.
\end{enumerate}
}
\end{prob}

%\ww
%{\bf Acknowledgment:}
%The author would like to express his sincere thanks to Izumi Ojima 
%for his interest in this topic and raising the above question.

%%%%%%%%%%%%%%%%%%%%%%%%%%%%%%%%%%%%%%%%%%%%%%%%%%%%%%%%%%%%%
%

\end{document}